\documentclass[12pt]{article}
\usepackage{latexsym}

\def\bs{\begin{subequations}}
\def\es{\end{subequations}}

\catcode`\@=11

\newtoks\@stequation
\def\subequations{\refstepcounter{equation}
  \edef\@savedequation{\the\c@equation}%
  \@stequation=\expandafter{\theequation}
  \edef\@savedtheequation{\the\@stequation}
  \edef\oldtheequation{\theequation}%
  \setcounter{equation}{0}%
  \def\theequation{\oldtheequation\alph{equation}}}

\def\endsubequations{\setcounter{equation}{\@savedequation}%
  \@stequation=\expandafter{\@savedtheequation}%
  \edef\theequation{\the\@stequation}\global\@ignoretrue}

\makeatletter
        \renewcommand{\theequation}{\thesection.\arabic{equation}}%
        \@addtoreset{equation}{section}%
\makeatother

\renewcommand{\thefootnote}{\fnsymbol{footnote}}

\begin{document}

\begin{titlepage}

 Revised November 12, 2012

\begin{center}

{\large \bf More Special Functions Trapped\\}

\vskip 0.2in

Charles Schwartz\footnote{E-mail: schwartz@physics.berkeley.edu}

{\em Department of Physics,
     University of California\\
     Berkeley, California 94720}
        
\end{center}

\vskip .3in



\vfill

\begin{abstract}
We extend the technique of using the Trapezoidal Rule for efficient 
evaluation of the
Special Functions of Mathematical Physics given by integral 
representations. This technique was recently used for 
Bessel functions, and here we treat Incomplete Gamma functions
and the general Confluent Hypergeometric Function.

\end{abstract}

\vfill

\end{titlepage}

\renewcommand{\thefootnote}{\arabic{footnote}}
\setcounter{footnote}{0}
\renewcommand{\thepage}{\arabic{page}}
\setcounter{page}{1}

\section{Introduction}
In a recent work I was able to demonstrate the power of the 
Trapezoidal Rule for numerical integration \cite{CS1},
\begin{equation}
\int_{-\infty}^{\infty}\; dx\; f(x) \approx \sum_{n}\;h\;f(nh),
\end{equation}
applied to integral 
representations of the various Bessel functions. \cite{CS2} Here, $h$ 
is the mesh size, which will be systematically reduced to show the 
convergence of the method; and the sum over mesh points $n$, 
nominally from $-\infty$ to $+\infty$, will be truncated when the 
contributions to the sum are below the desired accuracy. Now that 
general technique is applied to other members of the venerable family 
called the Special Functions of Mathematical Physics.

First, we look at the Incomplete Gamma function (the earlier paper 
did deal briefly with the complete Gamma function); and we find that 
a particular approach initiated some years ago by Talbot \cite{T}, 
employing the inverse Laplace transform, is very useful.

Then, following more traditional lines, we start with the textbook integral 
representation for the general Confluent Hypergeometric function and 
change variables so that it behaves remarkably well under the 
Trapezoidal Rule.

Lots of numerical results are presented, showing the rapid 
convergence and versatility of the calculational method.

\section{Incomplete Gamma Function}

The definiton is
\begin{equation}
\gamma(s,x) = \int_{0}^{x}\;dt\;t^{s-1}\;e^{-t};\label{a1}
\end{equation}
and I want to follow Talbot's\cite{T} idea, further developed by 
 Weideman and Trefethen \cite{WT} \cite{ST}, about looking at 
inverse Laplace Transforms.

Taking the Laplace transform of (\ref{a1}), we find
\begin{equation}
\int_{0}^{\infty}\;dx e^{-yx}\;\gamma(s,x) = 
\Gamma(s)\;y^{-1}\;(1+y)^{-s}.\label{a2}
\end{equation}
The inverse transform gives us back the original function:
\begin{equation}
\gamma(s,x) = \frac{\Gamma(s)}{2\pi i} \int dy \;
e^{xy}\;y^{-1}\;(1+y)^{-s} =  \frac{\Gamma(s)}{2\pi i} \int dy \;
e^{y}\;y^{-1}\;(1+y/x)^{-s},\label{a3}
\end{equation}
where we have assumed, for now, that $x$ is a positive real number, 
so that the second equality above is just the result of scaling the 
integration variable $y$.

The contour of this y-integral starts way out in the third quadrant 
 and ends far out in the second quadrant, crossing the real axis  
 at a positive 
value of y. This is the same contour we earlier used \cite{CS2} for calculating 
the inverse of the complete Gamma function. We find that for $x \rightarrow 0$ 
this gives us just that earlier formula; although the original 
definition, Eq. (\ref{a1}), gives the complete Gamma function when $x 
\rightarrow \infty$.

I have not seen this integral representation (\ref{a3}) for $\gamma(s,x)$ before. 
It should be great for 
using the Trapezoidal Rule. We start by 
choosing a simple contour: 
\begin{equation}
y =c + 1 - cosh(u) + i\; sinh(u),\;\;\;\;\; -\infty < u < \infty ;
\end{equation}
and $c$ is where the contour crosses the real axis. We see that the 
integral (\ref{a3}) has a simple pole at $y=0$ and a branch point at 
$y=-x$.
For $c > 0$ we have the function $\gamma(s,x)$; and if we take $-x < 
c < 0$, then we get the answer $ \gamma(s,x) - \Gamma(s) 
= -\int_{x}^{\infty}\;dt\;t^{s-1}\;e^{-t}$.

Looking for the point of constant phase to be at this crossing, we 
find the formula
\begin{equation}
c = [(s+1-x) + \sqrt{(s+1-x)^{2}+4x}\;]/2
\end{equation}
The data in Tables 1 and 2 below show some examples of the ratio 
$P(s,x)=\gamma(s,x)/\Gamma(s)$ computed by this method. These 
calculations, as in Ref\cite{CS2}, are done with standard double 
precision accuracy ($10^{-16}$) and I truncate the sum over mesh 
points when the fractional contributions are less than this amount.

\vskip 0.5cm
Table 1. Computations of $P(s,x)$  using the Trapezoidal Rule

\begin{tabular} {||c|c|c|c||}\hline
1/h & $P(0.1,1)\times 10$ & $P(1,0.1)\times 100 $ & $P(0.1,0.1)\times 
10 $ \\ \hline
1 & 9.85296 26362 19827 & 9.58023 11961 90712  & 8.37021 11048 74736   \\ \hline
2 & 9.75973 88897 94659 & 9.51192 39220 99238  & 8.27666 75413 58269   \\ \hline
4 & 9.75872 65150 00294 & 9.51625 80387 73782  & 8.27551 75852 50319   \\ \hline
8 & 9.75872 65627 36719 & 9.51625 81964 04048  & 8.27551 75958 58505   \\ \hline
16& 9.75872 65627 36723 & 9.51625 81964 04037  & 8.27551 75958 58505   \\ \hline
\end{tabular}
\vskip 0.5cm
\newpage 

Table 2. Computations of $P(s,x)$  using the Trapezoidal Rule

\begin{tabular} {||c|c|c|c||}\hline
1/h & $P(1,1)\times 10$ & $P(10,10)\times 10 $ & 
$P(1000,1000)\times 10 $ \\ \hline
1 & 6.35418 14003 89701 &5.69571 39376 33220    & 5.45885 48876 16142\\ \hline
2 & 6.32017 05027 55139 &5.41549 12109 30272    & 5.11064 21436 17747\\ \hline
4 & 6.32120 56442 73888 &5.42070 15780 49574     & 5.04658 65045 16275\\ \hline
8 & 6.32120 55882 85574 &5.42070 28552 84053     & 5.04204 40307 32861\\ \hline
16& 6.32120 55882 85577 &5.42070 28552 81477     & 5.04205 21812 85238\\ \hline
32&                     &5.42070 28552 81479     & 5.04205 24418 02230\\ \hline 
64&                       &                      & 5.04205 24418 02222 \\ \hline 
\end{tabular}
\vskip 0.5cm
The number of mesh points used for the last line of data in Table 1 
was 153,151,153 respectively; and for Table 2 the numbers were  
151, 285, and 841.

Special cases of the Incomplete Gamma Function include the Error 
Function and various Exponential Integrals.

\section{Confluent Hypergeometric Functions}

The focus here is on the integral representation,
\begin{equation}
C(a,b;x) =\int_{0}^{1}dt\;t^{a-1}\;(1-t)^{b-1}\; e^{xt}, \;\;\; a, b 
> 0.\label{c1}
\end{equation}
This is connected to the traditional definition of the Confluent 
Hypergeometric Function, as follows.
\begin{equation}
_{1}F_{1}(a,c=b+a;x) = M(a,c=b+a;x) = B(a,b)^{-1}\;C(a,b;x).
\end{equation}
And this introduces the classical Beta-function,
\begin{equation}
B(a,b) = \Gamma(a)\;\Gamma(b)/\Gamma(a+b) = C(a,b;x=0).
\end{equation}
where that last identity comes from the fact that the Hypergeometric 
functions are defined to have the value 1 at $x=0$.

We propose to evaluate these functions by the Trapezoidal Rule, after 
making some convenient (real) coordinate transformations under the 
integral:
\begin{equation}
t = \frac{1}{2}(1+tanh(u)) = e^{u}/(e^{u}+e^{-u}); \;\;\;\;\; u = 
sinh(v);\label{c2}
\end{equation}
and then we apply the Trapezoidal rule to the infinite integral over v.

Some numerical results are shown in Tables 3-7  below.

\vskip 0.5cm
Table 3. Computations of $C(a=1.,b=1.;x)$ 

\begin{tabular} {||c|c|c|c||}\hline
1/h & $(x=0.)$ & $(x=1.) $ & $(x=100.)\times 10^{-41} $ \\ \hline
1 &1.00107 86347 90329  &1.72266 65011 24113  &1.34247 60795 12472 \\ \hline
2 &1.00000 01397 11616  &1.71828 29887 33384  &2.83640 72153 82483 \\ \hline
4 &1.00000 00000 00001  &1.71828 18284 59079  &2.68831 33551 31161 \\ \hline
8 &                     &1.71828 18284 59044  &2.68811 71418 07843 \\ \hline
16&  &                                        &2.68811 71418 16129 \\ \hline
32&  &                                        &2.68811 71418 16129 \\ \hline 
\end{tabular}
\vskip 0.5cm

Table 4. Computations of $C(a=0.1,b=1.;x)$ 

\begin{tabular} {||c|c|c|c||}\hline
1/h & $(x=0.)$ & $(x=1.) $ & $(x=100.)\times 10^{-41} $ \\ \hline
1 & 9.99991 08876 77476  &11.21464 29773 5867  &1.34425 85914 11300 \\ \hline
2 &10.00000 00993 10149  &11.21300 57977 0100  &2.86472 18192 35488 \\ \hline
4 & 9.99999 99999 99998  &11.21300 52032 3319  &2.71295 94923 56232 \\ \hline
8 & 9.99999 99999 99998  &11.21300 52032 3318  &2.71278 37414 62587 \\ \hline
16&  &                                         &2.71278 37414 71210 \\ \hline
32&  &                                         &2.71278 37414 71210 \\ \hline 
\end{tabular}
\vskip 0.5cm

Table 5. Computations of $C(a=0.1,b=10.;x)$ 

\begin{tabular} {||c|c|c|c||}\hline
1/h & $(x=0.)$ & $(x=1.) $ & $(x=100.)\times 10^{-29} $ \\ \hline
1 &7.53951 55733 97154  & 7.63169 79919 69269 &3.40378 23912 96282 \\ \hline
2 &7.59131 58970 84419  & 7.67046 16609 55571 &1.70396 29450 35887 \\ \hline
4 &7.59138 00009 01282  & 7.67049 54154 30269 &1.09112 36709 31219 \\ \hline
8 &7.59138 00009 11013  & 7.67049 54154 32864 &1.07365 08998 41708 \\ \hline
16&7.59138 00009 11017  & 7.67049 54154 32878 &1.07365 07978 79342 \\ \hline
32&                     &                     &1.07365 07978 79343 \\ \hline 
\end{tabular}

\newpage 

Table 6. Computations of $C(a=10.,b=0.1;x)$ 

\begin{tabular} {||c|c|c|c||}\hline
1/h & $(x=0.)$ & $(x=1.) $ & $(x=100.)\times 10^{-44} $ \\ \hline
1 &7.53951 55733 97154  &20.26525 88483 9334  &1.69844 72080 59734 \\ \hline
2 &7.59131 58970 84419  &20.44048 93266 6513  &1.59374 46727 15385 \\ \hline
4 &7.59138 00009 01285  &20.44076 89724 1024  &1.59966 19996 26905 \\ \hline
8 &7.59138 00009 11014  &20.44076 89724 7923  &1.59966 08127 76379 \\ \hline
16&7.59138 00009 11021  &20.44076 89724 7924  &1.59966 08127 76238 \\ \hline
32&  &                     &1.59966 08127 76246 \\ \hline 
\end{tabular}

\vskip 0.5cm
Table 7. Computations of $C(a=0.1,b=0.1;x)$ 

\begin{tabular} {||c|c|c|c||}\hline
1/h & $(x=0.)$ & $(x=1.) $ & $(x=100.)\times 10^{-44} $ \\ \hline
1 &19.71352 25754 7283  &35.95446 58322 5812  &1.70487 63130 76804 \\ \hline
2 &19.71463 97119 7631  &35.95643 50355 3144  &1.61024 26167 51076 \\ \hline
4 &19.71463 94890 5016  &35.95643 47587 2009  &1.61504 43786 53350 \\ \hline
8 &19.71463 94890 5015  &35.95643 47587 2014  &1.61504 16242 89936 \\ \hline
16&                     &35.95643 47587 2013  &1.61504 16242 89858 \\ \hline
32&  &                                        &1.61504 16242 89859 \\ \hline 
\end{tabular}
\vskip 0.5cm
The number of mesh points needed to obtain the best results shown in 
the Tables above was always under 200.

Looking at this data, one is reminded of a particular virtue  of the 
present method.\cite{CS1} The integral in Eq. (\ref{c1}) has, for non-integral 
values of a and b, singularities at the end points $t=0,1$. This 
means that traditional methods of numerical integration (Simpson's 
rule, Richardson extrapolation, Gauss quadrature, etc.) would not 
work at all well. But the transformation Eq. (\ref{c2}) used for our 
approach with the Trapezoidal Rule, handles those singularities very 
nicely: they are moved to infinity and smothered beneath exponential 
decays.

The same method used above should be applicable to the general 
Hypergeometric function,
\begin{equation}
_{2}F_{1}(a,b;c;z) = B(b,c-b)^{-1}\; \int_{0}^{1}\;dt\; t^{b-1} \;
(1-t)^{c-b-1} \;(1-zt)^{-a},
\end{equation}
although I have not done those calculations. There are various 
transformation formulas for $_{2}F_{1}$ that let one deal with the 
singularities when $z$ approaches 1 or infinity.

\section{Discussion}

The numerical results displayed above show very nice rates of 
convergence using the Trapezoidal Rule. The programming is straightforward; 
and there is freedom 
for the user to explore various avenues.

While some other authors, e.g. \cite{WT}, have focused on finding 
optimal contours, my own opinion is that this general method is so 
powerful and robust that such refinements may be more a distraction than a benefit.

Since I have restricted myself here to real values of the coordinates 
(x) and the
parameters (s,a,b), it is left for others to explore the extension of 
those variables  into the complex planes.  

I am grateful to J.A.C.Weideman for some stimulating conversations.


\begin{thebibliography}{99}
	
\bibitem{CS1}
C. Schwartz, ``Numerical integration of analytic functions,'' J. 
Comp. Phys.  \textbf{4} (1969) 19. 

\bibitem{CS2}
C. Schwartz, ``Numerical Calculaton of Bessel Functions,'' Int. J. 
Mod. Phys. C \textbf{23} (12), 1250084 (2012); arXiv:1209.1547  

\bibitem{T}
A. Talbot, ``The accurate numerical inversion of Laplace transforms,''
J. Inst. Math. Appl. \textbf{23} (1979) 97.

\bibitem{WT}
J. A. C. Weideman and L. N. Trefethen, ``Parabolic and hyperbolic 
contours for computing the Bromwich integral,'' Mathematics of 
Computation \textbf{76} (2007) 1341.

\bibitem{ST}
T. Schmelzer and L. N. Trefethen, ``Computing the Gamma Function 
Using Contour Integrals and Rational Approximations,'' SIAM J. 
Numer.Anal. \textbf{45} (2007) 558.
 
\end{thebibliography}
\end{document}